\documentclass[11pt,twoside]{article}
\usepackage{latexsym}
\usepackage{amssymb,amsbsy,amsmath,amsfonts,amssymb,amscd,amsthm}
\usepackage{graphicx,color}
\usepackage{hyperref}

\def \longrightharpoonup {\relbar\joinrel\rightharpoonup}

\newcommand{\fer}[1]{(\ref{#1})}

\newcommand{\commentout}[1]{}
\newcommand{\R}{\mathbb{R}}

\newcommand {\al} {\alpha}

\newcommand {\e}  {\varepsilon}

\newcommand {\vp} {\varphi}

\newcommand {\Chi} {{\bf \raise 2pt \hbox{$\chi$}} }

\newcommand {\car} { {\mathcal R} }
\newcommand {\f}   {\frac}
\newcommand {\p}   {\partial}
\newcommand{\dis}{\displaystyle}
\newcommand{\beq}{\begin{equation}}
\newcommand{\eeq}{\end{equation}}
\newcommand{\bea} {\begin{array}{rl}}
\newcommand{\eea} {\end{array}}
\newcommand{\bc} {\begin{cases}}
\newcommand{\ec} {\end{cases}}
\newcommand{\bepa}{\left\{ \begin{array}{l}}
\newcommand{\eepa} {\end{array}\right.}
\newtheorem{theorem}{Theorem}[section]
\newtheorem{lemma}[theorem]{Lemma}

\title{\Large \bf Time fluctuations in a population model of adaptive dynamics}

\author{
Sepideh Mirrahimi\thanks{
CNRS, Institut de Math\'ematiques de Toulouse UMR 5219, 31062 Toulouse, France. Email: sepideh.mirrahimi@math.univ-toulouse.fr}
\thanks{ Universit\'e de Toulouse ; UPS, INSA,  UT1, UTM ; IMT ; 31062 Toulouse, France.}
\and  Beno\^ \i t Perthame\footnotemark\thanks{UPMC Universit\'e Paris 06, Laboratoire Jacques-Louis Lions, CNRS UMR7598, 4 pl. Jussieu, 75005 Paris, France}
\thanks{INRIA EPI BANG and Institut Universitaire de France . Email: benoit.perthame@upmc.fr}
\and Panagiotis E. Souganidis\thanks{The University of Chicago,
Department of Mathematics, 5734 S. University Avenue,
Chicago, IL 60637, USA. Email: souganidis@math.uchicago.edu}
\thanks{Partially supported by the National Science Foundation
}}

\date{\today}

\begin{document}
\maketitle
\pagestyle{plain}
\pagenumbering{arabic}

\begin{abstract}
\noindent We study the dynamics of  phenotypically structured populations in environments with fluctuations. In particular, using novel arguments from the theories of Hamilton-Jacobi equations with constraints and homogenization, 
we obtain results about the evolution of populations in environments with time oscillations, the development of concentrations in the form of Dirac masses, the location of the dominant traits and their evolution in time.  Such questions have already been studied in time homogeneous  environments.
More precisely we consider  the dynamics of a phenotypically structured population in a changing environment under mutations and competition for a single resource. The mathematical model is a non-local parabolic equation with a periodic in time reaction term.  
We study the asymptotic behavior of the solutions in the limit of small diffusion and fast reaction. Under concavity assumptions on the reaction term, we prove that the solution converges to a Dirac mass whose evolution in time is driven by a Hamilton-Jacobi equation with constraint and an effective growth/death rate which is derived as a homogenization limit. We also prove that, after long-time, the population concentrates on a trait where the maximum of an effective growth rate is attained. Finally we provide an example showing that the time oscillations may lead to a strict  increase of the asymptotic population size.
 
\end{abstract}

\noindent {\bf Key-words}: Reaction-diffusion equations, Asymptotic analysis, Hamilton-Jacobi equation, Adaptive dynamics, Population biology, Homogenization.
\\[2mm]
{\bf AMS Class. No}: 35B25, 35K57, 49L25, 92D15

\section{Introduction}
\label{sec:int}

Phenotypically structured populations can be modeled using non-local Lotka-Volterra equations, which  have the property that, in the small mutations limit,  the solutions concentrate on one or several evolving in time Dirac masses.
A recently developed mathematical approach, which  uses Hamilton-Jacobi equations with constraint,  allows us to understand the behavior of the solutions in constant environments \cite{OD.PJ.SM.BP:05,GB.BP:08,GB.SM.BP:09,AL.SM.BP:10}.  

\noindent Since stochastic and periodic modulations are  important for the modeling \cite{RL:76,SR:80,AS.SE:95,EK.RK.NB.SL:05,HS.CR.JH:11}, a natural and relevant question is whether it is possible to further develop the theory to models with time fluctuating environments. \\

\noindent In this note we consider  an environment which varies periodically in time in order, for instance, to take into account the effect of seasonal variations in the dynamics, and we study the asymptotic properties of the initial value problem 
\beq \begin{cases}
\e \, n_{\e,t}= n_\e  \; R \big(x, \f{t}{\e}, I_\e(t) \big) + \e^2 \Delta n_\e \quad \text{ in }  \quad \R^N \times ( 0,\infty),
\\[2mm]
n_{\e}(\cdot ,0)= n_{0,\e} \quad \text{ in } \quad \R^N, \\[2mm]
I_\e(t) := \int_{\R^N} \psi(x) n_\e(x,t) dx, 
\end{cases}
\label{model}
\eeq 
where 
\beq 
\label{asR}
R:\R^N \times \R \times [0, \infty) \to \R \quad  \text{is smooth and $1$-periodic in its second argument}.
\eeq

\noindent The population is structured by phenotypical traits $x\in \R^N$  with  density $n_\e(x,t)$ at time $t$. It is  assumed that there exists a single type of resource which is consumed by each individual trait $x$ at a rate $\psi(x)$;  $I_\e(t)$ is then the total consumption of the population. The mutations and the growth rate are represented respectively by the  Laplacian  term and $R$. The novelty is the periodic in time dependence of the growth rate $R$.
The small coefficient  $\e$ is used  to consider only rare mutations and to rescale time in order to study a time scale much larger than the  generation one.\\

\noindent {To ensure the survival and the boundedness of the population we assume that $R$ takes positive values for ``small enough populations'' and negative values for ``large enough populations" ,  i.e.,  there exist  a value $I_M>0$ 
such that 
\beq
\max_{0\leq s \leq 1, \;x\in \R^N} R(x, s, I_M) =0  \quad\text{and }Ê\quad  {\mathcal X} := \{x\in \R^N, \; \int_0^1 R(x,s,0)ds >0 \} \neq \emptyset.
\label{as2}
\eeq
 }

\noindent In addition the growth rate $R$ satisfies, for some positive constants $K_i$, $i=1,\ldots,7,$ and all $(x,s,I)\in \R^N\times \R\times [0,I_M]$ and  $A>0$,  the following concavity and decay assumptions:
\beq
-K_1 \leq D_x^2 R(x,s,I) \leq - K_2,  \qquad   K_3- K_1|x|^2 \leq R(x,s,I) \leq  K_4 - K_2|x|^2, 
\label{as1}
\eeq
\beq
-K_5 \leq D_I R(x, s, I) \leq -K_6,
\label{as11}
\eeq
\beq 
D_{x}^3 R \in L^\infty \big( \R^N \times (0,1)\times [0,A] \big) \quad \text{ and } \quad |D^2_{x,I}R|\leq K_7.
\label{as1b}
\eeq

\noindent The ``uptake coefficient''   $\psi :\R^N\to \R$ must be regular and bounded from above and below, i.e.,  there exist positive constants $\psi_m$, $\psi_M$ and $K_8$ such that
\beq
0 < \psi_m \leq \psi \leq \psi_M\quad\text{and}\quad \|\psi\|_{C^2} \leq K_8.
\label{as3}
\eeq

\noindent We also assume  that the initial datum is  ``asymptotically monomorphic'', i.e., it is close to a Dirac mass in the sense that there exist $x^0\in {\mathcal X}$,  $\rho^0>0$ and 
a smooth $u_\e^0:\R^N \to \R$ such that 
\beq
n_{0,\e}= e^{u_\e^0/ \e}  \  \   \text{and, \ as }  \e\to0,  
\eeq
\beq
n_\e(\cdot,0)  \underset{\e \to 0}{\longrightarrow} \varrho^0 \delta(\cdot-x^0) \quad \text{weakly in the  sense of measures}.
\label{as3b}
\eeq

\noindent In addition there exist constants  $L_i >0, i=1, \ldots,4,$  and a smooth $u^0:\R^N\to \R$ such that, for all $x\in\R^N$,
\beq
- L_1 I\leq  D^2_x u_\e^0\leq - L_2 I,\quad - L_3- L_1|x|^2\leq u_\e^0(x)\leq  L_4- L_2|x|^2, \quad  \max_{x\in \R^N} u^0(x) =0=u^0( x^0) 
\label{as4}
\eeq
and, as $\e\to0$,
\beq
u_\e^0 
{\longrightarrow}  
u^0 \quad \text{locally uniformly in \ $\R^N$}.
\label{as5}
\eeq

\noindent Finally, it is necessary to impose the following compatibility assumption between the initial data and the growth rate $R$:
\beq
\label{as:com}
4L_2^2\leq K_2\leq K_1\leq 4L_1^2.
\eeq

\noindent Our first result is about the behavior of the $n_\e$'s as $\e \to 0$. It asserts the existence of a fittest trait $\overline x(t)$ and a total population size $\overline \rho (t)$ at time $t$ and provides a ``canonical equation'' for the evolution in time of $\overline x$ in terms of the ``effective fitness'' ${\mathcal R}(x,y)$ 
satisfying  
${\mathcal R}(x,x)=0.$ In the sequel, $D_1{\mathcal R}$ denotes the derivative of ${\mathcal R}$ with respect to the first argument.

\begin{theorem}[Limit as $\e\to 0$] 
Assume \eqref{asR}--\eqref{as:com}.  There exist a fittest trait $\overline x\in \mathrm{C}^1\left([0,\infty); {\mathcal X} \right)$ and a total population size $\overline \rho\in \mathrm{C}^1\left([0,\infty);(0,\infty)\right)$ such that,  along subsequences  $\e\to 0$,    
$$  
n_\e(\cdot,t) {\longrightharpoonup}  \,\overline  \varrho(t) \delta(\cdot - \overline x(t)) \quad \text{weakly in the sense of measures,}
$$
and, 
$$
I_\e 
{\longrightharpoonup}\,  \overline  I := \overline  \varrho  \psi( \overline x)  \quad  \text{ in \ $L^\infty(0,\infty)$ weak-$\star$}.
$$
Moreover, $\overline x$ satisfies the canonical equation
\beq 
\label{canonical2}
\dot{\overline x}(t) =\left(-D^2_x u(\overline x(t),t) \right)^{-1} \cdot D_1 {\mathcal R}(\overline x(t),\overline x(t)).
\eeq
\label{th1}
\end{theorem}

\noindent We note that, in the language of adaptive dynamics,  $\mathcal R(y,x)$ can be interpreted as the effective fitness of a mutant $y$ in a resident population with a dominant trait $x$, while 
$D_1 \mathcal R$ is usually called the selection gradient, since  it  represents the capability of invasion. The extra term $\left(-D^2_x u(\overline x(t),t) \right)^{-1}$ is an indicator of the diversity around the dominant trait in the resident population.\\

\noindent The second issue is the identification of the long time limit of the fittest trait $\overline x$ .  We prove that, in the limit $t \to \infty$, the population converges to 
a, so called, Evolutionary Stable Distribution (ESD) corresponding to a distribution of population which is stable under introduction of small mutations (see \cite{JM.GP:73,IE:83,PJ.GR:09} for a more detailed definition). See also \cite{LD.PJ.SM.GR:08,Raoul2009-3} for recent studies of the local and global stability of stationary solutions of integro-differential population models in constant environments.
\begin{theorem}[Limit as $t\to \infty$] In addition to \eqref{asR}--\eqref{as:com} assume that either $N=1$ or, if $N>1$, $R$ is given, for some smooth $b, d, B, D:\R^N\to (0,\infty)$ by  
\beq
R(x,s,I) = b(x) B(s,I)-d(x) D(s,I). 
\label{specialR}
\eeq
Then, as $t\to \infty$, the population reaches an Evolutionary Stable Distribution $ \overline  \varrho_\infty \delta(\cdot - \overline x_\infty)$, i.e.,  
$
\overline  \varrho(t)  \,{\longrightarrow} \, \overline  \varrho_\infty \  \text{and} \    \overline x(t) \,  {\longrightarrow} \,\overline x_\infty,
$
where  $\overline  \varrho_\infty>0$ and $\overline x_\infty$ are characterized by ($\mathcal I$ is defined in \fer{periodic})
\beq
\label{maxRR}
{\mathcal R}( \overline x_\infty,\overline x_\infty) = 0 = \max_{x\in \R^N} {\mathcal R}( x, \overline x_\infty) \ \ \text{ and} \ \   \overline  \varrho_\infty=\f{1}{\psi( \overline x_\infty)}\int_0^1 \mathcal I( \overline x_\infty,s)ds.
\eeq 
\label{th:longtime}
\end{theorem}

\noindent Notice that we do not claim the uniqueness of the Evolutionary Stable Distribution. Indeed there may exist several $(\rho_\infty,\overline x_\infty)$ satisfying \fer{maxRR}.  Here we only prove that there exists $(\rho_\infty,\overline x_\infty)$ satisfying \fer{maxRR} such that,  as $t\to \infty$, the  population converges to $ \overline  \varrho_\infty \delta(\cdot - \overline x_\infty)$.\\

\noindent The difference between our conclusions and  the results for time homogeneous  environments in \cite{AL.SM.BP:10} is that, in the canonical equation \fer{canonical2}, the growth rate $R$ is replaced by an effective growth rate $\mathcal R$ which is derived after a homogenization process. Moreover, we are only able to prove that the  $I_\e$'s  converge in  $L^\infty$ weak-$\ast$ and not a.e., which  is  the case for constant environments in \cite{AL.SM.BP:10}. This adds a difficulty in Theorem \ref{th:longtime} and it is the reason why we are not able to describe, without additional assumptions, the long-time limit behavior of the fittest trait $\overline x$ for general growth rate $R$ when $N>1$. This remains an open question.\\

\noindent In Section \ref{sec:counter}, we give an example of $\mathcal R$ not satisfying the assumptions of Theorem \ref{th:longtime} for which $\overline x(t)$ exhibits a periodic behavior. This example fits the structure \fer{HJ} below with general concavity properties on $\mathcal R$ but it is not necessarily derived from a homogenization limit.\\

\noindent The proofs use in a fundamental way the classical Hopf-Cole transformation 
\beq\label{HC}
 u_\e = \e \, \ln \,n_\e, 
\eeq
 which yields the  following Hamilton-Jacobi equation for $u_\e$  :
\beq
\begin{cases} u_{\e,t} =  R \big(x, \f{t}{\e}, I_\e(t) \big) + |D_x u_\e|^2 + \e \Delta u_\e & \text{in $ \ \R^N\times(0,\infty),$}\\[2mm]
u_\e(\cdot,0)=u_\e^0,& \text{in \ $\R^N$.}
\end{cases}
\label{hje}
\eeq

\noindent The next theorem describes the behavior of the $u_\e$'s, as $\e\to 0$ (recall that $\mathcal R$ is defined in Section \ref{sec:proof}).
\begin{theorem} 
Assume \fer{asR}--\fer{as:com}. Along subsequences  $\e\to 0$, $u_\e\to u$ locally uniformly in $\R^N\times [0,\infty)$, where $u\in \mathrm{C}(\R^N\times [0,\infty))$ is a solution of
\beq\label{HJ}
\begin{cases}
 u_t = {\mathcal R}(x, \overline x(t) ) + | D_x u|^2 & \text{in \ $\R^N\times(0,\infty)$},
\\[2mm]
\dis \max_{x \in \R^N} u(x,t)= 0 = u(\overline x(t),t) & \text{in  \ $(0,\infty)$},
\\[2mm]
u(\cdot,0)=u^0 & \text{in \ $\R^N$.}
\end{cases}
\eeq
\label{th:hj}
\end{theorem}

\noindent In general not much is  known about the structure of the effective growth rate ${\mathcal R}$. 
In Section \ref{sec:ncon} we give an example of an $R$, for which the effective growth rate $\mathcal R$ can be computed explicitly. Moreover, for this example no concavity assumption is made.\\

\noindent Note that the convergence of the $u_\e$'s in Theorem \ref{th:hj} and, thus, the convergence of the $n_\e$'s in Theorem \ref{th1} are established only along subsequences. To prove convergence for all $\e$, we need that \fer{HJ} has a unique solution.  This is, however, not known even 
for non oscillatory  environments except for some
particular form of growth rate $R$ (see \cite{GB.BP:07,GB.BP:08}).\\

\noindent The paper is organized as follows. In Section \ref{sec:proof} we study the asymptotic behavior of the solution under rare mutations (limit $\e\to 0$) and we provide the proofs of Theorems \ref{th1} and \ref{th:hj}. In Section \ref{sec:ltb} we consider  the long time behavior of the dynamics (limit  $t\to \infty$) and  we give  the proof of Theorem \ref{th:longtime}. In Section \ref{sec:ncon} we study a particular form of growth rate $R$, for which the results can be proved without any concavity assumptions on $R$ and the effective growth rate $\mathcal R$ has a natural structure. Finally we present an example of an oscillatory environment which yields an asymptotic effective population density that is strictly larger than the averaged one.

\section{The behavior as $\e\to 0$ and the proofs of Theorems \ref{th1} and \ref{th:hj}}
\label{sec:proof}

We present the proofs of Theorems \ref{th1} and \ref{th:hj} which are closely related. Since the argument is long, we next summarize briefly the several steps. First we obtain some a priori bounds on $u_\e$ and $I_\e$. Then we identify the equation for the fittest trait $\overline x$. The limit of the $I_\e$'s is studied in Lemma \ref{lm:periodic}. The last three steps are the identification (and properties) of the effective growth rate $\mathcal R$, the effective Hamilton-Jacobi equation and the canonical equation. 


\proof[Proofs of Theorems \ref{th1} and \ref{th:hj} ]

\noindent{\em \large Step 1: a priori bounds.} It follows from \eqref{as2}, arguing as in  \cite{AL.SM.BP:10}, that for all $t\geq 0$,
\beq
0 < I_\e(t) \leq I_M+ O(\e).
\label{est:I}
\eeq 
Next we use \eqref{as1} and \fer{as4} to get, for some  $C_1,\,C_2>0$, $\e\leq 1$ and all $(x,t)\in \R^N\times [0,\infty)$,
\beq
- L_1  \leq  D_x^2u_\e  \leq - L_2\quad \text{and} \quad -L_3-L_1|x|^2-C_1t\leq u_\e(x,t)\leq L_4-L_2|x|^4+C_2t. 
\label{est:D2u}
\eeq
It follows from \fer{hje} that,  for all balls $B_R$ centered at the origin and of radius $R$, there exists $C_3=C_3(R)>0$ such that
\beq
 \| u_{\e,t} (x,t) \|_{ L^\infty(B_R \times (0,\infty))} \leq C_3.
\label{est:dtu}
\eeq
Finally the regularity properties of the ``viscous'' Hamilton-Jacobi equations yield that, for all $T>0$, there exists $C_4=C_4(R,T)>0$ such that
\beq
\| D_x^3 u_\e \|_{L^\infty (B_R\times [0,T])} \leq C_4.
\label{est:d3u}
\eeq
Hence, after differentiating \fer{hje} in $x$, the previous estimates yield a $C_5=C_5(R,T)>0$ such that
\beq
\|D^2_{t,x} u_\e\|_{L^\infty (B_R\times [0,T])} \leq C_5.
\label{est:dxtu}
\eeq
All the above bounds allow us to pass to the limit, along subsequences $\e \to 0$,
and to obtain $u: \R^N\times [0,\infty)\to \R$ such that, as $\e\to 0$,
$$
 u_\e \,{\longrightarrow}\, u \  \  \text {in \  \  $C_{\rm loc}(\R^N\times [0,\infty))$  \  \  and, for all $T>0$ and  $(x,t)\in \R^N \times [0,T],$ }
$$
$$
- L_1 \leq  D_{x}^2u(x,t)  \leq - L_2 \ \ \text{ and } \ \ u, D^2_{t,x} u,  D_{x}^3 u \in L^\infty (\R^N\times [0,T]).
$$

\noindent{{\em \large Step 2. The fittest trait.} In view of the strict concavity of $u_\e$, for each $\e>0$, there exists a unique $\overline x_\e(t)$ such that
$$
u_\e( \overline x_\e(t),t) = \max_{x\in \R^N} u_\e(x,t)  \  \  \text{ and } \  \  D_x u_\e( \overline x_\e(t),t)=0.
$$
Differentiating the latter equality with respect to $t$ and using \fer{hje} we find
$$
\dot{\overline  x}_\e(t) \cdot D_x^2 u_\e( \overline x_\e(t),t) = - D_x u_{\e,t}( \overline x_\e(t),t)
$$
$$
=- D_x R\big(\overline x_\e(t), \f{t}{\e}, I_\e(t)\big) - 2 D_x^2 u_\e( \overline x_\e(t),t) \cdot D_x u_\e( \overline x_\e(t),t) -\e \Delta D_x u_\e( \overline x_\e(t),t)
$$
$$
=- D_x R\big(\overline x_\e(t), \f{t}{\e}, I_\e(t)\big) -\e \Delta D_x u_\e( \overline x_\e(t),t).
$$

\noindent Since $D_x^2 u_\e( \overline x_\e(t),t)$ is invertible and $\| D_x^3 u \|_{L^\infty(B_R\times [0,T] )}\leq C_4$, it follows that $\dot {\overline  x}_\e (t)$ is bounded in $(0,T)$, and, hence,  along subsequences $\e\to 0$, \ 
${\overline  x_\e } {\rightarrow} {\overline  x }  \  \text{in }  C_{\rm loc} ((0,\infty)), $   for some   \  $\overline  x  \in C^{0,1}((0,\infty)) $ \   such that 
\beq\begin{cases}
u( \overline x(t),t)=\max_{x\in \R^N} u(x,t) , \qquad  D_x u( \overline x(t),t)=0, 
\\[2pt]
\text {and}\\[2mm]  
\dot {\overline  x }(t) = \big( -D^2_x u({\overline  x }(t),t)    \big)^{-1}  \cdot D_x \Big \langle R\big(\overline x(t), \f{t}{\e}, I_\e(t)\big) \Big \rangle,
\end{cases}
\label{canonical}
\eeq
where the bracket denotes the weak limit of $R\big(\overline x(t), \f{t}{\e}, I_\e(t)\big)$ which exists, since $R$ is locally bounded.
\\
\smallskip

\noindent{{\em \large Step 3. The weak limit of $I_\e$.} 
To identify the weak limit of the $I_\e$'s,  we consider the first exit time  $T^*>0$  of  $\overline x$  from  ${\mathcal X}$,  i.e., the smallest time $T^*>0$ such that $\overline x (t) \in {\mathcal X}$ for all $0 \leq t < T^*$ and $ \overline x(T^*) \in \p {\mathcal X}$ if $T^* <\infty$. Note that $T^*$ is well defined since $\overline x(0)=x^0 \in {\mathcal X}$. The last step of the ongoing proof is to show  that $T^*= \infty$.
\\

\noindent We need  the following two  results. Their proofs are given after the end of the ongoing one.  
\begin{lemma}
Assume \fer{as2}, \fer{as11} and \fer{as3}. For all $x\in {\mathcal X}$, there exists a unique $1$-periodic positive solution $ {\mathcal I} ( x,s ):[0,1] \to(0,I_M) $ to  
\beq \begin{cases}
\f{d}{d s} {\mathcal I} ( x,s ) =  {\mathcal I} ( x,s )  \; R \big(x,s ,  {\mathcal I} ( x,s ) \big),\\[2mm]
 {\mathcal I} ( x,0 )=  {\mathcal I} ( x, 1).
\end{cases}
\label{periodic}
\eeq
Moreover, as  ${\mathcal X} \ni x \to x_0 \in \p {\mathcal X}$, 
\beq
\max_{0\leq s \leq 1}  {\mathcal I} ( x,s )  \to 0 .
\label{frX}
\eeq
\label{lm:periodic}
\end{lemma}

\noindent In view of \fer{frX}, for $x \in \p {\mathcal X}$, we define, by continuity, $ {\mathcal I} ( x,s ) =0$.

\begin{lemma} Assume \fer{as2}, \fer{as11} and \fer{as3}. {Let $T_\e^*$  be the smallest time $T_\e^*>0$ such that $\overline x_\e (t) \in {\mathcal X}$ for all $0 \leq t < T_\e^*$ and $ \overline x(T^*_\e) \in \p {\mathcal X}$ if $T^*_\e <\infty$. Then, for all $0<t<T_\e^\ast$,
$$  \Big|  \ln I_\e(t) - \ln {\mathcal I} ( \overline x_\e (t) , \f t \e  ) \Big|  \leq  \Big|  \ln I_\e(0) - \ln {\mathcal I} ( \overline x_\e (0) , 0 ) \Big|e^{-\f{K_6 t}{\e}}+ C {\sqrt \e } ,$$
where $C$ only depends on the constants $K_i$. Moreover, as $\e\to 0$, $T_\e^*\to T^*$.}
Consequently, if $ \overline x(t) \in {\mathcal X}$ for $0\leq t <T^*$ and $ \overline x(T^*) \in \p {\mathcal X}$, then, as $\e \to 0$  and $ t \to T^* $, $I_\e(t)  {\longrightarrow} 0.$
%
\label{lm:periodic3}
\end{lemma}
\noindent It follows that, as $\e\to 0$,  
\beq 
I_\e(\cdot)  {\longrightharpoonup}  \overline I(\cdot) =  \int_0^1{\mathcal I} ( \overline x (\cdot) ,s) ds >0 \qquad  \text{ in $L^\infty((0,T^*))$ weak-$\star$} .
\label{weaklI}
\eeq
Once $\overline I$ is known, it is possible to compute the weight $\overline \varrho $ of the Dirac mass. Indeed, we show in the next steps that, as $\e\to 0$,
$$
I_\e(\cdot)=\int_{\R^N} \psi(x) n_\e(x,\cdot) dx \,{\longrightharpoonup}\,  \overline  I(\cdot) = \overline \varrho (\cdot) \psi(\overline x (\cdot)) \ \  \text{ in $L^\infty((0,T^*))$ weak-$\star$}. 
$$

\smallskip

\noindent{{\em \large Step 4. The effective growth rate.} 
We can now explain the average used to determine the effective growth rate. Again \eqref{as11} and Lemma \ref{lm:periodic3} yield that, as $\e \to 0$,
$$
\int_0^{T^*}  \Big| R \big(x,\f t \e ,  I_\e ( t ) \big) - R \big(x,\f t \e ,  {\mathcal I} ( \overline x (t) , \f t \e )  \big) \Big| dt\leq K_5 \int_0^{T^*} \big|  I_\e ( t )  -  {\mathcal I} ( \overline x (t) , \f t \e  )  \big| dt   {\longrightarrow} 0.
$$
Therefore the weak limit in \eqref{canonical} is computed as the weak limit (in time)  of $R \big(x,\f t \e ,  {\mathcal I} ( \overline x (t) , \f t \e )  \big)$. To this end, we define, for all $x \in R^N \ \text{and} \ y\in \overline {\mathcal X}$ (here ${\overline {\mathcal X}}$ stands for the closure of ${\mathcal X}$),
\beq \label{def:cR2}
{\mathcal R}( x,y):=\int _0^1R(x,s,{\mathcal I} ( y , s ) )ds.
\eeq
It follows that, for $0\leq t \leq T^*$, 
\beq \label{def:cR}
\Big \langle R\big( x, \f{t}{\e}, I_\e(t)\big) \Big \rangle = {\mathcal R}( x, \overline x(t))=\int_0^1 R \big(x,s ,  {\mathcal I} ( \overline x (t) , s )  \big) ds.
\eeq
In particular, if $y \in \p {\mathcal X}$, then  ${\mathcal R}( x,y) = \int_0^1 R(x,s,0) ds$. 
\\

\noindent Notice also that, integrating  \eqref{periodicJ} below in $s$ and using the  periodicity, we always  have, for $x\in \overline {\mathcal X}$,
\beq 
{\mathcal R}( x, x ) \equiv 0.
\label{Rvanishes}
\eeq
Finally,  it is immediate from  \fer{as1} and \fer{def:cR2},  that $\mathcal R(x,y)$ is strictly concave in the first variable.
\\

\noindent{{\em \large Step 5. The limiting Hamilton-Jacobi equation.} 
It is now possible  to pass to the limit $\e\to 0$ in \eqref{hje} for $(x,t) \in \R^N \times [0, T^*)$. To this end, observe that 
$$\vp_\e(x,t):=u_\e(x,t)-\int_0^t R\left(x,\f \tau \e,I_\e(\tau)\right)d\tau$$
solves 
$$\vp_{\e,t}-\e\Delta \vp_\e=\e \int_0^t \Delta R\left(x,\f \tau \e,I_\e(\tau)\right)d\tau+|D_x \vp_\e+\int_0^t D_x R\left(x,\f \tau \e,I_\e(\tau)\right)d\tau|^2.$$
Since the $u_\e$'s  converges locally uniformly from Step 1 and $R\left(x,\f t \e,I_\e(t)\right)$ converges weakly in $t$ and strongly in $x$ to $\mathcal R\left(x,\overline x(t)\right)$ from Step 2, we find that, as $\e \to 0$,
$$\vp_\e(x,t){\longrightarrow} \vp(x,t)=u(x,t)-\int_0^t \mathcal R\left(x,\overline x(\tau)\right)d\tau \quad \text{in \ $C_{\rm loc}(\R^N \times (0, T^*))$}.$$
Moreover, in view of  \fer{as1},  \fer{as1b} and Lemma \ref{lm:periodic}, for any $T \in (0,T^*)$ and $R>0$, there exists $C=C(R,T)>0$  such that, for $(x,t) \in  B_R \times [0,T]$
$$
\left|\int_0^t \Delta R\left(x,\f \tau \e,I_\e(\tau)\right)d\tau \right|\leq C,
$$
and, as $\e \to 0$,
$$
\int_0^t |D_x R\left(x,\f \tau \e,I_\e(\tau)\right)d\tau-D_x R\left(x,\f \tau \e,\mathcal I(\overline x(\tau),\f{\tau}{\e})\right)|d\tau\leq K_7 \int_0^T \big|  I_\e ( t )  -  {\mathcal I} ( \overline x (t) , \f t \e  )  \big|dt  \;  {\longrightarrow} \; 0.
$$
Thus, as $\e\to 0$ and for all $(x,t) \in B_R \times [0,T]$,
$$
\int_0^t D_x R\left(x,\f \tau \e,I_\e(\tau)\right)d\tau{\longrightarrow} \int_0^t D_1\mathcal R\left(x,\overline x(\tau)\right)d\tau.
$$

\noindent It follows from the stability of viscosity solutions that $\vp$ is a viscosity solution to
$$ \vp_t =\big | D_x \vp+ \int_0^t D_1 \mathcal R\left(x,\overline x(\tau)\right)d\tau\big |^2 \  \text{in } \  \R^N\times (0,T^*),$$
which, written in terms of  $u$, reads
$$ u_t =\mathcal R\big( x , \overline x(t)\big) + | D_x u|^2  \  \text{in \ $\R^N\times (0,T^*)$} .$$
The constraint $\max_{x\in \R^N}u(x,t)=0$ follows from \fer{HC} and \fer{est:I} (see \cite{GB.BP:08,GB.SM.BP:09}). We then conclude following \cite{GB.BP:08} that, for all $t\in (0,T^\ast)$,
$$
n_\e(\cdot,t) {\longrightharpoonup}  \,\overline  \varrho(t) \delta(\cdot - \overline x(t)) \quad \text{weakly in the sense of measures.}
$$
}
\\

\noindent{{\em \large Step 6. The canonical equation.} 
The canonical equation \eqref{canonical2} now follows from  \eqref{canonical} and  \eqref{def:cR}.
\\

\noindent{{\em \large Step 7. The global time $T^*=\infty$.} 
Assume $T^*<\infty$. Then $ \overline x(T^*) \in \p {\mathcal X}$.  It follows from the canonical equation \eqref{canonical2} that, for all $t \in (0,T^*)$, 
$$\bea
\f d{dt} \int_0^1 R(  \overline x(t) ,s, 0) ds & = \dot{  \overline x} (t) \int_0^1 D_xR(  \overline x(t) ,s, 0) ds
\\[5pt]
& = D_1\mathcal R\big( {\overline x(t)} , \overline x(t)\big)  \big(-D^2u(\overline x(t),t)\big)^{-1}  \int_0^1 D_xR(  \overline x(t) ,s, 0) ds,
\eea$$
while, when $t=T^*$, Lemma \ref{lm:periodic3} yields that $D_1\mathcal R\big(  {\overline x(t)} , \overline x(t)\big)= \int_0^1  D_x R(  \overline x(t) ,s, 0) ds$. 
Hence
$$
\f d{dt} \int_0^1 R(  \overline x(T^*) ,s, 0) ds >0,
$$
which  is a contradicition because, by the definition of the open set ${\mathcal X}$,  $\int_0^1 R(  \overline x(t) ,s, 0) ds>0$ for $t \in [0,T^*)$ and 
$\int_0^1 R(  \overline x(T^*) ,s, 0) ds =0$.
\qed

\bigskip

\proof[Proof of Lemma \ref{lm:periodic}]  First  we prove that, for a fixed $x \in {\mathcal X}$, there exists a solution $\mathcal I$ of \fer{periodic}. To this end observe that  ${\mathcal J}: = \ln {\mathcal I}$ solves 
\beq \begin{cases}
 \f{d}{ds}{\mathcal J} ( x,s ) =   R \big(x,s ,  \exp({\mathcal J} ( x,s ) \big)   \  \text{in \  $s\in[0,1]$},
\\[2mm]
{\mathcal J} ( x,0 )=  \al. 
\end{cases}
\label{periodicJ}
\eeq
It turns out that it is possible to choose  $\al \leq  \ln I_M$ so that ${\mathcal J} ( x,0 )={\mathcal J} ( x,1 )$.
Indeed, the definition of $I_M$ in  \eqref{as2} yields that, if 
$\al= \ln I_M$, then ${\mathcal J}(x,1) < \al .$  On the other hand,
for $\al$ very small we claim that ${\mathcal J} (x,1) > \al$, which is enough  to conclude, since 
${\mathcal J} ( x,s)$ been  a continuous increasing function of $\al$, it has a fixed point $\al^\ast$. Choosing $\al=\al^\ast$ yields a periodic solution.\\

\noindent To prove the claim, we set $\mu = \int_0^1 R(x,s, 0) ds >0$ since $x \in {\mathcal X}$. Because $R$ is locally bounded, there exists  a constant $C>0$, which is   independent of $\al$,  such that  ${\mathcal J}(x,s) \leq \al +C$ and,  for $\al$ small enough,
$$
{\mathcal J}(x,1) = {\mathcal J}(x,0) + \int_0^1 R(x, s,  \exp({\mathcal J} ( x,s ) \big)   ds \geq {\mathcal J}(x,0) +  \int_0^1 R(x, s,  0 \big)   ds + O(e^{\al}) \geq  {\mathcal J}(x,0) + \f{\mu}{2}.
$$
This proves the claim and the existence of a periodic solution.
\\

\noindent The uniqueness follows from a contraction argument. Indeed let $\mathcal J_1$ and $\mathcal J_2$ be two periodic solutions to \fer{periodicJ}. Then
$$
\f{d}{d s} (\mathcal J_1-\mathcal J_2)=R \big(x,s ,  \exp({\mathcal J_1} ( x,s ) \big) -R \big(x,s ,  \exp({\mathcal J_2} ( x,s ) \big) .
$$
Multiplying the above equation by $\mathrm{sgn}\;(\mathcal J_1-\mathcal J_2)$ and  using the monotonicity in $I$ according to  \fer{as11}, we find
$$
\f{d}{d s} \big |\mathcal J_1-\mathcal J_2\big |\leq -C |\mathcal J_1-\mathcal J_2\big |,
$$
and, after integration,
$$
C \int_0^1 |\mathcal J_1(s)-\mathcal J_2(s)\big |ds\leq -|\mathcal J_1(1)-\mathcal J_2(1)|  + |\mathcal J_1(0)-\mathcal J_2(0)| = 0,
$$
and, hence, $\mathcal J_1=\mathcal J_2$.
\\

\noindent Finally we prove \eqref{frX}. It follows from \eqref{periodicJ} that, for $x \in {\mathcal X}$,  
$$
0= \int_0^1 R(x, s, e^{{\mathcal J}(x,s)} ) ds \leq  \int_0^1 R(x, s , 0) ds - K_6 e^{\min_{0\leq s \leq 1} {\mathcal J}(x,s)}. 
$$
If  $x \to x_0 \in \p {\mathcal X}$, then $ \int_0^1 R(x, s , 0) ds \to 0$ and, since  the variations of ${\mathcal J}(x,s)$ are bounded, because  $R$ is locally bounded, the result follows. 
\qed

\proof[Proof of Lemma \ref{lm:periodic3}]  
We identify the weak limit of $I_\e$ and prove \fer{weaklI}. We begin with the observation that in the ``gaussian''- type concentration,  $x- \overline x_\e(t)$ scales as $\sqrt \e$. 
\\

\noindent Indeed multiplying  \eqref{model} by $\psi$ and integrating with respect to $x$ we find (recall that with $J_\e :=\ln I_\e$),
$$
\e \f{d}{d t} J_\e ( t ) = \f{\int_{\R^N} \psi(x) n_\e(x,t)   \; R \big(x,\f t \e ,  I_\e ( t )  \big) dx}{\int_{\R^N} \psi(x) n_\e(x,t)  dx}  
+ \e^2 \f{ \int_{\R^N} \Delta\psi (x)  \; n_\e(x,t) dx}{\int_{\R^N} \psi(x) n_\e(x,t)  }.
$$
Note that in order to justify the integration by parts above, we first replace $\psi$ by $\psi_L=\chi_L \psi$ where  $\chi_L$ is a compactly supported smooth function such that $\chi_L \equiv 1$ in ${\rm{B}(0,L)}$ and  $\chi_L \equiv 0$ in $\R^N\backslash\rm{B}(0,2L)$. Then we integrate by parts and finally let $L\to +\infty$. 

\noindent Returning to the above equation we find 
$$
\e \f{d}{d t} J_\e ( t ) = \f{\int_{\R^N} \psi(x)e^\f{u_\e(x,t)-u_\e(\overline x_\e(t),t)}{\e}  \; R \big(x,\f t \e ,  I_\e ( t )  \big) dx}{\int_{\R^N} \psi(x) e^\f{u_\e(x,t)-u_\e(\overline x_\e(t),t)}{\e}  dx}  + O(\e^2)
$$
$$
= \f{\int_{\R^N} \psi(x)e^\f{u_\e(x,t)-u_\e(\overline x_\e(t),t)}{\e}  \; \big[R \big(x,\f t \e ,  I_\e ( t )  \big) - R \big( \overline x_\e(t) ,\f t \e ,  I_\e ( t )  \big) \big] dx}{\int_{\R^N} \psi(x) e^\f{u_\e(x,t)-u_\e(\overline x_\e(t),t)}{\e}  dx}  +R \big( \overline x_\e(t) ,\f t \e ,  I_\e ( t )  \big)   + O(\e^2).
$$

\noindent Using Laplace's method for approximation of integrals, \eqref{est:D2u}  and \fer{est:d3u}, we find that  the first term is of order $\sqrt \e$ and, hence, 
$$
\e \f{d}{d t} J_\e ( t ) =R \big( \overline x_\e(t) ,\f t \e ,  I_\e ( t )  \big)   + O(\sqrt \e).
$$

\noindent Next we compute
\begin{equation*}\begin{split}
\e \f{d}{d t} [{\mathcal J} ( \overline x_\e(t),  \f t \e ) - J_\e(t) ]&= R \Big( \overline x_\e(t) ,\f t \e ,  \exp \big({\mathcal J} ( \overline x_\e(t),  \f t \e )\big) \Big)  - R \big( \overline x_\e(t) ,\f t \e ,  \exp( J_\e ( t ))  \big)  \\
&+ O(\sqrt \e) + \e \,D_x {\mathcal J} ( \overline x_\e(t),  \f t \e )\,\dot {\overline  x}_\e (t).
\end{split}
\end{equation*}
\noindent Multiplying the above equality  by ${\rm sgn}({\mathcal J} ( \overline x_\e(t),  \f t \e ) - J_\e(t))$, using our previous estimates  and employing the monotonicity property in \eqref{as11},  we get
$$
\e \f{d}{d t} \big| {\mathcal J} ( \overline x_\e(t),  \f t \e ) - J_\e(t) \big|= -\big|  R \big( \overline x_\e(t) ,\f t \e ,  \exp({\mathcal J} ( \overline x_\e(t),  \f t \e ) ) \big)   - R \big( \overline x_\e(t) ,\f t \e ,  \exp( J_\e ( t ) ) \big) \big|  + O(\sqrt \e)
$$
$$
\leq - K_6 \big| {\mathcal J} ( \overline x_\e(t),  \f t \e ) - J_\e(t) \big|  + O(\sqrt \e). 
$$
The first claim of Lemma \ref{lm:periodic3} is now immediate. {Moreover, since $\overline x_\e(t)\to \overline x(t)$, locally uniformly as $\e\to 0$, we obtain that, as $\e \to  0$,
$$T_\e^*\to T^*.$$
The last claim is a consequence of the previous steps and Lemma \ref{lm:periodic}.
}
\qed

\section{The long time behavior}
\label{sec:ltb}

\subsection{Convergence as $t \to \infty$ when $N=1$  (The proof of Theorem \ref{th:longtime} (i))}

Throughout this subsection we assume that $N=1$. The goal is to prove the existence of some $\overline x_\infty \in \R^N$ such that, as $t\to \infty$, $\overline x(t)\rightarrow \overline x_\infty$ and
\beq
\label{ess}
\car(\overline x_\infty,\overline x_\infty)=0=\max_{x\in R} \car(x,\overline x_\infty).
\eeq
To this end, we consider the map 
$A: \R\longrightarrow \R$  defined by  $A(x)=y$, where $y$ is the unique maximum point of $\car(\cdot,x)$.
We obviously have
$$D_x \car(A(x),x)=0.$$

\noindent We consider the following three cases depending on the comparison between  $\overline x(\cdot)$ and $A(\overline x(\cdot))$.
If  $\overline x(t)<A(\overline x(t))$, then $D_x \car(\overline x(t), \overline x(t))>0$ and, if $\overline x(t)>A(\overline x(t))$, then $D_x \car(\overline x(t), \overline x(t))<0$. It then follows using  \fer{canonical2} and the concavity of $u$ that, if $\overline x(t) < A(\overline x(t))$ (resp. $\overline x(t) > A(\overline x(t))$), then $\dot{\overline x}(t)>0$ (resp. $\dot{\overline x}(t)<0$).
 If  $\overline x(t)=A(\overline x(t))$, then   $D_x \car(\overline x(t), \overline x(t))=0$
and hence again \fer{canonical2} yields $\dot{\overline x}(t)=0$.
We also notice that, if $A(\overline x(0))=\overline x(0)=x_\infty$, then from the above argument we have that  for all $t\geq 0$, $\overline x(t)=x_\infty$, with $x_\infty$ satisfying \fer{ess}.
\smallskip

\noindent Now we assume that $A(\overline x(0))>\overline x(0)$ (the case $A(\overline x(0))<\overline x(0)$ can be treated similarly) and set 
$$t_0:=\inf\, \left\{t\in \R \; : \; A(\overline x(t)) \,\leq\, \overline x(t) \right\}.$$
If $t_0<\infty$, then $A(\overline x(t_0))=\overline x(t_0)$ and, hence, for all $t\geq t_0$, $\overline x(t)=\overline x(t_0)=x_\infty$.
If $t_0=\infty$, then $\dot{\overline x}(t)>0$ for all $t\geq 0$, and thus, since the set $B=\{\overline x(t)\;:t\;\in [0,\infty)\}$ is compact (see below), there exists $\overline x_0\in \R$ such that
$$\lim_{t\to \infty}\overline x(t)=\overline x_0.$$
The compactness of $B$ follows from the observation that, in view of \fer{def:cR2}, \fer{as1} and \fer{as11}, 
$$
\mathcal{R}(x,\overline x(t))\leq K_4-K_2|x|^2,
$$
and, since $\mathcal{R}(\overline x(t),\overline x(t))=0$, 
\beq\label{xbound}|\overline x(t)|\leq (K_4/K_2)^{1/2}.\eeq
 We now claim that $\overline x_0$ satisfies \fer{ess}. Indeed, if there exists $z\in \R^N$ such that $\car(z,\overline x_0)>0$, then using \fer{HJ}, we have
$\lim_{t\to \infty }u(z,t)=+\infty,$
a contradiction to the constraint $\max_{x\in \R} u(x,t)=0$.

\subsection{Convergence for a particular case  with $N>1$ (The proof of Theorem \ref{th:longtime} (ii))}

We prove Theorem \ref{th:longtime} in the multi-d case with a growth rate $R$ as  in \fer{specialR}.  In this case we find
\beq
\label{md}
 b(\overline x(t))  \langle B(s,I(t)  \rangle  - d(\overline x(t))   \langle  D(s,I(t)  \rangle =0,
\eeq
and thus 
$$
\big( -D^2_x u({\overline  x }(t),t)    \big) \dot{\overline x}(t) =  b'({\overline  x }(t))  \langle B(s,I(t)  \rangle  - d'({\overline  x }(t))   \langle  D(s,I(t)  \rangle
$$
$$
=  \left[ \f{ b'({\overline  x }(t)) }{ b({\overline  x }(t))} - \f{ d'({\overline  x }(t)) }{ d({\overline  x }(t))}  \right]  d({\overline  x }(t))  \langle  D(s,I(t)  \rangle.
$$
Therefore,  after taking inner product with $ \dot{\overline x}(t)$, dividing by $d({\overline  x }(t))  \langle  D(s,I(t)  \rangle$ and using the strict concavity of $u$, we obtain
$$
| \dot{\overline x}(t)|^2 \leq C \f{d}{dt} \ln \left( \f{b({\overline  x }(t))}  {d({\overline  x }(t))} \right).
$$
\noindent This proves that $t\mapsto b({\overline  x }(t))/ d({\overline  x }(t)) $ increases and thus converges, as $t \to \infty$,  to some constant $l$. For this we need to show that $\{\overline x(t):t\in [0,\infty)\}$ is bounded, a fact which follows exactly as in the proof of \fer{xbound}.  Finally, in view of \fer{md}, we also have
$$\lim_{t\to\infty} \f{ \langle  D(s,I(t)  \rangle }{ \langle B(s,I(t)  \rangle }=l.$$

\noindent We prove next that
$$l=\max_x \f{b(x)}{d(x)}.$$
Arguing by contradiction we assume that $l<\max_x (b(x)/d(x))$. Then there must exist  $\widetilde x\in \R^N$ such that
$l<(b(\widetilde x)/d(\widetilde x))$, in which case
$$
 0<\liminf_{t\to \infty}b(\widetilde  x)  \langle B(s,I(t)  \rangle  - d(\widetilde x)   \langle  D(s,I(t)  \rangle.
$$
Finally, since $u$ solves 
$$
\p_t u =|D_x u|^2 + b(x) \langle B(s,I(t)  \rangle  - d(x)   \langle  D(s,I(t)  \rangle,
$$
we find 
$$\lim_{t\to \infty}u(\widetilde x,t)=\infty,$$
a contradiction to the constraint $\max_{x\in \R^N} u(x,t)=0$.\\

\subsection{A counterexample in the multi-dimension case}\label{sec:counter}

In this subsection we present an example  showing that, when $N>1$,  the $\overline x$ 's  may not converge, as $t \to \infty$,  at least for the Hamilton-Jacobi problem \fer{HJ}. Indeed we find a strictly concave with respect to the first variable ${\mathcal R}: \R^N \times \R^N \to \R$, an $1$-periodic map $t \to \overline x(t)$ and 
a function 
$u:\R^N \times [0,\infty]\to \R$ which satisfies 
$$
\begin{cases}
\p_t u-|D_x u|^2=\mathcal{R}\left(x,\overline x(t)\right) \ \ \text{ in } \ \  R^N \times (0,\infty), \\[2mm]
\max_{x\in \R^N}u(x,t)=u\left(\overline x(t),t\right)=0,\\[2mm]
u(\cdot, 0)=u_0 \ \text{ in } \  \R^N.
\end{cases}\label{equc}
$$

\noindent We choose $G : \R^N \to \R^N$ so that the ode $\dot{x}=G(x)$ has a periodic solution; note that such function exists only for $N>1$. A simple example for $N=2$ is
$G(x_1,x_2)=(-x_2,x_1),$
which admits  $(x_1(t),x_2(t))=(r\cos t,r\sin t)$ as periodic solutions.\\

\noindent Let $F:\R^N\to \R$ be an arbitrary smooth function and define $\mathcal{R}:\R^N\times \R^N\to \R$ by
\beq\label{CR}\mathcal{R}(x,y)=-\left(D F(y)G(y)+4F(y)^2\right)|x-y|^2+2F(y)G(y)(x-y).
\eeq
It is immediate that $\mathcal R$ is a concave function with respect to $x$ and satisfies $\mathcal R(x,x)=0$. It is also easily verified that
$$u(x,t)=-F(\overline x(t))|x-\overline x(t)|^2\qquad \text{with}\quad \dot{\overline x}(t)=G(\overline x(t)) \; \text{and} \; \overline x(0)=x_0,$$
is a viscosity solution of \fer{equc} for $\mathcal{R}$ 	as in \fer{CR} and 
$u_0(x)=-F(x_0)|x-x_0|^2.$
Moreover the canonical equation \fer{canonical} is written as
$$\dot{\overline x}(t)=\left(-D^2_x u(\overline x(t),t)\right)^{-1}D_x \mathcal{R}(\overline x(t),\overline x(t))=(2F(\overline x(t))^{-1}\left(2F(\overline x(t))G(\overline x(t))\right)=G(\overline x(t)).$$
Finally we choose $x_0 \in \R^N$ such that $t\mapsto \overline x(t)$ with $x(t_0)=x_0$ is $1$-periodic. Then  the limit $\lim_{t\to \infty} \overline x(t)$ does not exist.\\

\noindent Note that the counterexample presented above is for the Hamilton-Jacobi problem \fer{HJ}. We do not  know if such periodic oscillation can arise in the  $\e \to 0$ limit of the viscous Hamilton-Jacobi equation \fer{hje}. When the  growth rate independent of time, a result similar to  Theorem \ref{th:longtime}  was proved in \cite{AL.SM.BP:10}  for general $R$. In that problem, the key point leading to the convergence, as $t \to \infty$, of the $\overline x(t)$'s is that $\overline I(t)$, which is the strong limit of the $I_\e$'s as $\e\to 0$,  is increasing in time. In the case at hand,  we can only prove that the $I_\e$'s  converge weakly to $\overline I$. We know nothing about the monotonicity of  $\overline I$. We remark that numerical computations suggest (see Figure \ref{fig:osc}) that monotonicity holds, if at all, in the average. 
\begin{figure}[h]
\begin{center}\includegraphics[angle=0,scale=0.4]{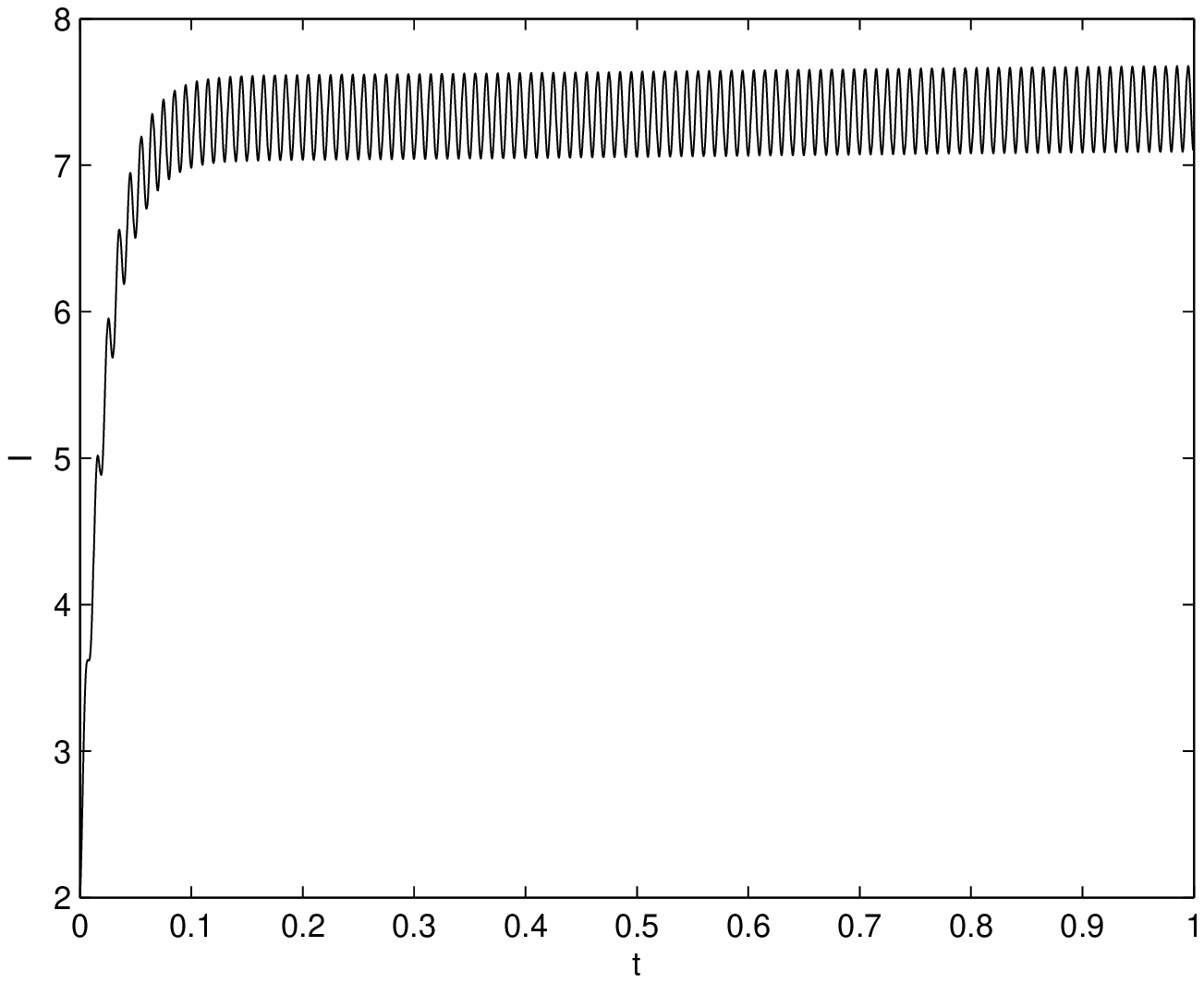}
\end{center}
\vspace{-25pt}
\caption{Dynamics of the total population $I_\e(t)$ for $R(x,s,I)=\left( 2+\sin\left(2\pi s \right) \right)\f{2-x^2 }{I+.5}-.5$, $\psi(x) =1$ and $\e=0.01$. The $I_\e$'s  oscillate with period of order $\e$ around a monotone  curve $\overline I$.}
\label{fig:osc}
\end{figure}

\section{A particular case with a natural structure for $\mathcal R$}
\label{sec:ncon}

The concavity assumption \fer{as1} is very strong. Here we study, using a different method based on BV estimates,  a class of  growth rates $R$ which do not satisfy \fer{as1}. Throughout this section, for several arguments, we follow \cite{GB.SM.BP:09} which studies a similar problem but without time oscillations.\\

\noindent We consider growth rates of the form
\beq\label{Rpart}
R(x,s,I)=b(x)B(s,I)-D(s,I),
\eeq
with 
\beq\label{Rpart0}
B, D : \R\times [0,\infty) \to R  \ \  \text{ $1$-periodic with respect to the first argument}
\eeq
and we assume  that, for all $(s,I)\in \R\times \big[\,\f{\widetilde I_m}{2},2\widetilde I_M\,\big]$ and $x\in \R^N$,
\beq\label{as:S}
0< B(s,I), \; 0<D(s,I)  \ \ \text{ and } \ \ 0<b_m\leq b(x)\leq b_M,
\eeq 
where $\widetilde I_M \, > \, \widetilde I_m \, > 0$ are such that
\beq \label{maxR}\max_{0\leq s \leq 1, \;x\in \R^N} R(x, s,\widetilde I_M) =0  \ \ \text{ and } \ \   \min_{0\leq s \leq 1, \;x\in \R^N} R(x, s,\widetilde I_m ) = 0, \eeq
and there constants $a_1>0$ and $a_2>0$ such that, for all $(s,I)\in \R\times [0,\infty)$,
\beq\label{as:SI}
\ D_I\,  B(s,I) \,<  -a_1  \ \text{ and} \    a_2 < \, D_I \, D(s,I).
\eeq

\noindent As far as $n_\e(\cdot,0)$ is concerned, we replace \fer{as3b}--\fer{as5} by
\beq\label{as:inibis}
\begin{array}{c}
\widetilde I_m\leq \int \psi(x) n_\e(x,0)dx\leq \widetilde I_M,  \text{ and}  \\[3mm]
 n_\e(x,0)\leq \exp\left(\f{-A|x|+B}{\e}\right)  \ \ \text{for some   $A,B>0$ and all $x\in \R^N$}.
\end{array}
\eeq

\noindent Eventhough \fer{Rpart} seems close to \fer{specialR}, no concavity assumption is made and the analysis of Section \ref{sec:ltb} does not apply here.\\
 
\begin{theorem}\label{th:sp}
Assume \fer{as3} and \fer{as:SI}--\fer{as:inibis}. Along subsequences $\e \to 0$,  the $u_\e$'s converge locally uniformly to $u \in \rm{C}(\R^N\times \R)$ satisfying the constrained Hamilton-Jacobi equation
\beq\label{HJ2}\begin{cases}
 u_t = {\mathcal R}(x,F(t) ) + | D_x u|^2 \ \ \text{ in} \ \ \R^N\times (0,\infty),
\\[2mm]
\dis \max_{x \in \R^N} u(x,t)= 0,\\[2mm]
\dis u(\cdot,0)=u^0 \ \text{ in } \ \R^N,
\end{cases}
\eeq
with 
$$
\mathcal R(x,F)=\int_0^1 \mathcal I(F,s)ds\left(\f{b(x)}{F}-1\right) \ \  \text{ and } \ \ F(t)=\lim_{\e \to 0}\, \dfrac{\int \psi(x) \, b(x) \, n_\e(x,t) \, dx}{I_\e(t)},
$$
and $\mathcal I$ defined in \fer{periodic2} below. In particular, along subsequences $\e \to 0$ and in the sense of measures,  $n_\e \, {\longrightharpoonup}\, n$
with $\mathrm supp\; n\subset \{(x,t) \,: \, u(x,t)=0\}\subset\{(x,t) \,:\, \mathcal R(x,F(t))=0\}$.
\end{theorem}

\noindent As in Theorem \ref{th:longtime}, we can deduce the long time convergence to the Evolutionary Stable Distribution. To this end we assume that

\beq\label{as:psimax}\text{there exists a unique  }x_\ast\in \R^N \text{ such that } b(x_\ast)=\max_{x\in\R^N}\;b(x).\eeq

\begin{theorem}\label{thm:longb}
Assume \fer{as3}, \fer{as:SI}--\fer{as:inibis} and \fer{as:psimax}
Then, as $t\to \infty$, the population reaches the Evolutionary Stable Distribution $\rho_\ast\delta(\cdot -x_\ast)$, i.e.,
\beq \label{nlimb}n(\cdot,t)\underset{t\to \infty}{\longrightharpoonup} \rho_\ast\delta(\cdot-x_\ast)  \ \  \text{in the sense of measures,}\eeq
with
$$\rho_\ast=\f{1}{\psi(x_\ast)}\int_0^1 \mathcal I\left(b(x_\ast),s\right)ds.$$
\end{theorem}
\proof[Proof of Theorem \ref{th:sp}]
It follows easily from \fer{maxR}, \fer{as:inibis} and the arguments in \cite{GB.SM.BP:09} that
 \beq \label{III} \widetilde I_m+O(\e)\leq I_\e(t)\leq \widetilde I_M+O(\e).\eeq
Define next
$$F_\e(t):=\f{\int b(x)\,\psi(x)\,n_\e(x,t)\,dx}{I_\e(t)},$$
and note that
\beq\label{boundF}b_m\leq F_\e(t)\leq b_M. \eeq

\noindent We next prove that $F_\e \in \text {BV}_{\text{loc}}(0,\infty)$ uniformly in $\e$. Indeed, using \fer{as3}, \fer{as:S} and \fer{III}, we find
\beq\label{fe}\begin{array}{rl}\f{d}{dt}F_\e(t)&={I_\e^{-2}}\left(I_\e\,\int n_{\e,t} \, b \,\psi \,dx-\int n_{\e,t} \, \psi\, dx\;\int n_\e\, b\, \psi \,dx\right)\\[2.5mm]
  &={I_\e^{-2}} ( I_\e\,\int \left(\e\Delta n_\e+ \e^{-1}{n_\e}\left(b \,B(\f t \e,I_\e)-D(\f t \e, I_\e)\right)\right) \, b\,\psi\,dx\\[2.5mm]
&-\int n_\e\, b\, \psi \,dx\, \int \big(\e\Delta n_\e+ \e^{-1}{n_\e}\left(b \,B(\f t \e,I_\e)-D(\f t \e,I_\e)\right) \psi\,dx)\\[2.5mm]
&=O(\e)+ (\e I_\e^2)^{-1}B(\f t \e,I_\e) \left(\int n_\e\, \psi\,dx\cdot \int n_\e\,b^2\, \psi\,dx -(\int n_\e \, b\, \psi\,dx)^2\right)\geq O(\e). \\
  \end{array}
\eeq
Then \fer{boundF} and \fer{fe}  yield that, for each $T>0$,  there exists  $C=C(T)>0$ such  that
$$\int_0^T |\f{d}{dt} F_\e|dt\leq C.$$
It follows that, along subsequences $\e\to 0$ , the $F_\e$'s converge a.e. and in $L^1$ to some $F$.

\noindent To conclude we need a result similar to the one  of  Lemma \ref{lm:periodic}
\begin{lemma}
Assume  \fer{as3} and \fer{as:SI}--\fer{as:inibis}. For all $t \in \R$, there exists a unique, $1$-periodic solution $  {\mathcal I} ( t,\cdot ) \in C^1(\R \to [\widetilde I_m, \widetilde I_M]$ to  
\beq \begin{cases}
\f{d}{d s} {\mathcal I} ( F( t ) , s ) =  {\mathcal I} \left( F(t) , s \right)  \; \big( F(t) B \left( s , {\mathcal I} ( F(t),s ) \right)- D \left(s, {\mathcal I} ( F(t) , s ) \right) \big), 
\\[2mm]
 {\mathcal I} ( F(t),0 )=  {\mathcal I} ( F(t),1). 
\end{cases}
\label{periodic2}
\eeq
Moreover,  for  all $T>0$ and as $\e\to 0$, \qquad
$
\int_0^T  |  I_\e(t) - {\mathcal I} ( F(t) , \f t \e  )| dt  \to  0. 
$
\label{lm:periodic2}
\end{lemma}

\noindent The first claim is proved as in Lemma \ref{lm:periodic}. We postpone the proof of the second assertion  to the end of this section.\\

\noindent Using \fer{as3}, \fer{as:SI}--\fer{as:inibis} and following \cite{GB.SM.BP:09}, we show that the $u_\e$'s  are  bounded and locally Lipschitz continuous  uniformly  in $\e$ and, hence, converge along subsequences  $\e \to 0$ to a solution $u$ of 
$$u_t=|D_x u|^2+ \int_0^1 B (s,  {\mathcal I} ( F(t) , s )) ds  b(x)-\int_0^1  D\left(s,{\mathcal I} (F(t),s) \right)ds.$$
Since $ {\mathcal I} (F, \cdot )$ is a periodic solution to \fer{periodic2}, we have
$$\int_0^1  B(s,  {\mathcal I} ( F(t) , s ) ds \  F(t)-\int_0^1   D (s,{\mathcal I} (F(t),s)) ds=0.$$
It follows that
$$ u_t=|D_x u|^2+\mathcal R(x,F(t))$$
with 
$$\mathcal R(x,F(t))=\int_0^1  D(s,{\mathcal I} (F(t),s) ) ds \ (\f{b(x)}{F(t)}-1).$$
The last claim of Theorem \ref{th:sp} can be proved using \fer{HC}, \fer{HJ2} and following \cite{GB.BP:08}.
\qed
\\

\noindent We conclude with
\proof[Proof of Theorem \ref{thm:longb}]
It follows from \fer{fe} that $F$ is an increasing function. Hence, in view of  \fer{boundF},  there exists $F_\ast$ such that,  as $t \to \infty$, \
$
F(t) \longrightarrow F_\ast.
$ \ 
Moreover
\beq
\label{Fast}
F_\ast=\max_{x\in \R^N}\; b(x)=b(x_\ast).
\eeq
Indeed, if not, then $\mathcal R(x_\ast,F_\ast)>0$ and, hence, from \fer{HJ2}, \ 
$
\lim_{t\to \infty} u( x_\ast,t)=\infty,
$
a contradiction to the constraint  $\max_{x\in \R^N} \;u(x,t)=0$. \\
Finally \fer{nlimb} follows from \fer{Fast} and the observation that 
 $$\mathrm supp\; n(x,t)\subset \{(x,t)\, :\,  \mathcal R(x,t)=0\}=\{(x,t) \,:\, b(x)=F(t)\}.$$
\qed

\proof[Proof of the second claim of Lemma \ref{lm:periodic2}]
Eventhough we follow the same ideas  as in  Lemma \ref{lm:periodic},  to prove the second claim, we need to modify the arguments, since, without the concavity assumption \fer{as1}, $u_\e$ may have several maxima.\\

\noindent  We use again the log transformations $J_\e =\log I_\e$ and $\mathcal J=\log \mathcal I$. Multiplying \eqref{model}  by $\psi(x)$ and integrating with respect to $x$ leads to
$$
\e \f{d}{d t} I_\e ( t ) =  B(\f t \e ,  I_\e ( t )) \int_{\R^N} \psi(x) n_\e(x,t)   \; b(x) dx- I_\e(t) D(\f t \e ,  I_\e ( t )) + \e^2  \int_{\R^N} \Delta\psi (x)  \; n_\e(x,t) dx.
$$
It follows from  \fer{periodic2} that
$$
\begin{array}{rl}
\e \f{d}{d t} \left( J_\e ( t ) - \mathcal J \left(F_\e(t), \f t \e \right)\right)&=F_\e(t) B\left(\f t \e ,  I_\e ( t )  \right)-F(t) B\left(\f t \e ,  \mathcal I \left(F_\e( t),\f t\e \right)  \right)
\\
\\
& - D\left(\f t \e ,  I_\e ( t )  \right) +D\left(\f t \e ,  \mathcal I \left(F_\e( t),\f t\e \right)  \right)+ O(\e^2).
\end{array}
$$
Multiplying the above equality  by $\rm{sgn}( J_\e(t) - {\mathcal J} ( F_\e(t),  \f t \e ) )$ and employing  \eqref{as:SI},  we obtain
$$
\begin{array}{rl}
\e \f{d}{d t} \left| J_\e ( t ) - \mathcal J \left(F_\e(t), \f t \e \right)\right|&=-\left| F_\e(t) B\left(\f t \e ,  I_\e ( t )  \right)-F_\e(t) B\left(\f t \e ,  \mathcal I \left(F_\e( t),\f t\e \right)  \right)\right|
\\
\\
& -\left|  D\left(\f t \e ,  I_\e ( t )  \right) -D\left(\f t \e ,  \mathcal I \left(F_\e( t),\f t\e \right)  \right) \right |+ O(\e^2).
\end{array}
$$
Integrating in time over $[0,T]$, for some fixed $T>0$, and using the convergence of the $F_\e$'s we find that, as $\e \to 0$,
$$
\int_0^T  | F(t) B(\f t \e ,  I_\e ( t ) ) - F(t) B(\f t \e ,  \mathcal I (F( t),\f t\e )) |+ 
| D(\f t \e ,  I_\e ( t ) )  - D(\f t \e ,  \mathcal I \big(F( t),\f t\e))|  dt\longrightarrow 0. 
$$
The second claim of  Lemma \ref{lm:periodic2} follows in view of \fer{as:SI}.
\qed

\section{A qualitative effect: fluctuations may increase the population size} 

 We conclude with an example that shows that the time-oscillations may lead to a strict increase of the population size at the evolutionary stable state, a conclusion which also holds in the context of physiologically structured populations \cite{CGP}. \\
 
\noindent  To this end, we consider, along the lines of Section \ref{sec:ncon}, the rate function 
\beq \label{example}
R(x,I)=b(x) - D_1(s)D_2(I)
\eeq
with $b$ and $D(s,I)=D_1(s)D_2(I)$ satisfying \eqref{Rpart0}--\eqref{as:SI} and \eqref{as:psimax} and, for simplicity, we take $\psi \equiv1$ in  \eqref{model}. The  goal is to compare the size of the ESD  in Theorem~\ref{thm:longb} to the one  obtained from the model with  the  ``averaged rate''
$$
R_{\text {av}}(x,I) = b(x) - D_{1,\text {av}}  D_2(I)  \  \text{ with} \   D_{1,\text {av}} = \int_0^1 D_1(s) ds. 
$$
Later, we write $f_{\text {av}}$ for the average of the  $1$-periodic map $f:\R\to \R$, i.e., $f_{\text {av}}=\int_0^1 f(s) ds$.
\\

\noindent Let  ${\mathcal I}$ be the $1$-periodic solution of \eqref{periodic2} with $F(t) \equiv b(x_\star)$ according to \eqref{Fast}.
With the above simplifications, the  magnitude $\rho_\star$ of the Evolutionary Stable Distribution obtained in \fer{thm:longb} is 
$$
\rho_\star=\int_0^1 {\mathcal I}(s) ds.
$$

\noindent Since we can multiply equation \eqref{periodic2} by any function of ${\mathcal I}(s)$, elementary maipulations  lead to the identities 
\beq \label{average1}
b(x_\star) = \int_0^1 D_1(s) D_2({\mathcal I}(s))ds \ \  \text{ and }  \ \    \int_0^1 D_2({\mathcal I}(s)) ds \;  b(x_\star) = \int_0^1 D_1(s) D_2^2 ({\mathcal I}(s)) ds .
\eeq
A straightforward application of the Cauchy-Schwarz inequality in \eqref{average1} yields
$$
b(x_\star)^2 \leq D_{1,\text {av}} \; D_2({\mathcal I})_\text{av} b(x_\star)
$$
and thus 
\beq \label{average 2}
b(x_\star) \leq D_{1,\text{av}} \int_0^1 D_2({\mathcal I}(s)) ds.
\eeq

\noindent Consider next the ``averaged'' version of \eqref{model}, i.e., the equation 
\beq  \begin{cases}
\e  n_{\e,\text{av},t}= n_{\e,\text{av}}  \; R_{\text{av}} (x, I_{\e,\text{av}}(t)) + \e^2 \Delta n_{\e,\text{av}} \quad \text{ in }  \quad \R^N \times ( 0,\infty),
\\[5pt]
 I_{\e,\text{av}}(t) =\int_{\R^N} n_{\e,\text{av}}(x,t) dx.
\end{cases}
\label{average_model} \eeq
It follows from the earlier results \cite{AL.SM.BP:10} that the magnitude $\rho_{\text{av}}$ of the Evolutionary Stable Distribution corresponding to \eqref{average_model}
satisfies the identity
\beq \label{average 3}
b(x_\star) = D_{1,\text{av}} D_2 (\rho_{\text{av}}),
\eeq
and, therefore, unless $D_1$  is  constant, in which case  \eqref{average 2} must be  an equality, we conclude
\beq \label{average 4}
D_2 (\rho_{\text{av}}) < \int_0^1 D_2({\mathcal I}(s) ds.
\eeq
\noindent If, in addition to above hypotheses, we also assume that 
\beq \label{D}
I \to D_2 (I) \ \  \text{ is concave},
\eeq
\noindent then \eqref{average 4} yields 
$$ \rho_{\text{av}} < \rho_\star,$$
which substantiates our claim about the possible effect of the time oscillations.


\begin{thebibliography}{10}

\bibitem{GB.SM.BP:09}
G.~Barles, S.~Mirrahimi, and B.~Perthame.
\newblock Concentration in {L}otka-{V}olterra parabolic or integral equations:
  a general convergence result.
\newblock {\em Methods Appl. Anal.}, 16(3):321--340, 2009.

\bibitem{GB.BP:07}
G.~Barles and B.~Perthame.
\newblock Concentrations and constrained {H}amilton-{J}acobi equations arising
  in adaptive dynamics.
\newblock {\em Contemp. Math.}, 439:57--68, 2007.

\bibitem{CGP}
J.~Clairambault, S.~Gaubert, and B.~Perthame.
\newblock Comparison of the perron and floquet eigenvalues in monotone
  differential systems and age structured equations.
\newblock {\em C. R. Acad. Sc. Paris}, 345((10):549--555, 2007.

\bibitem{LD.PJ.SM.GR:08}
L.~Desvillettes, P.-E. Jabin, S.~Mischler, and G.~Raoul.
\newblock On mutation-selection dynamics for continuous structured populations.
\newblock {\em Commun. Math. Sci.}, 6(3):729--747, 2008.

\bibitem{OD.PJ.SM.BP:05}
O.~Diekmann, P.-E. Jabin, S.~Mischler, and B.~Perthame.
\newblock The dynamics of adaptation: an illuminating example and a
  {H}amilton-{J}acobi approach.
\newblock {\em Th. Pop. Biol.}, 67(4):257--271, 2005.

\bibitem{IE:83}
Ilan Eshel.
\newblock Evolutionary and continuous stability.
\newblock {\em Journal of Theoretical Biology}, 103(1):99--111, 1983.

\bibitem{PJ.GR:09}
P.-E. Jabin and G.~Raoul.
\newblock On selection dynamics for competitive interactions.
\newblock {\em J. Math. Biol.}, 63(3):493--517, 2011.

\bibitem{EK.RK.NB.SL:05}
E.~Kussel, R.~Kishony, N.~Q. Balaban, and S.~Leibler.
\newblock Bacterial persistence: a model of survival in changing environments.
\newblock {\em Genetics}, 169:1807--1814, 2005.

\bibitem{RL:76}
R.~Lande.
\newblock Natural selection and random genetic drift in phenotypic evolution.
\newblock {\em Evolution}, 30:314--334, 1976.

\bibitem{AL.SM.BP:10}
A.~Lorz, S.~Mirrahimi, and B.~Perthame.
\newblock Dirac mass dynamics in multidimensional nonlocal parabolic equations.
\newblock {\em Comm. Partial Differential Equations}, 36(6):1071--1098, 2011.

\bibitem{JM.GP:73}
J.~Maynard~Smith and G.~R. Price.
\newblock The logic of animal conflict.
\newblock {\em Nature}, 246:15--18, 1973.

\bibitem{GB.BP:08}
B.~Perthame and G.~Barles.
\newblock Dirac concentrations in {L}otka-{V}olterra parabolic {PDE}s.
\newblock {\em Indiana Univ. Math. J.}, 57(7):3275--3301, 2008.

\bibitem{Raoul2009-3}
G.~Raoul.
\newblock Long time evolution of populations under selection and vanishing
  mutations.
\newblock {\em Acta Applicandae Mathematica}, 114:1--14, 2011.

\bibitem{SR:80}
S.~Rosenblat.
\newblock Population models in a periodically fluctuating environment.
\newblock {\em Journal of Mathematical Biology}, 9:23--36, 1980.

\bibitem{AS.SE:95}
A.~Sasaki and S.~Ellner.
\newblock The evolutionarily stable phenotype distribution in a random
  environment.
\newblock {\em Evolution}, 49(2):337--350, 1995.

\bibitem{HS.CR.JH:11}
H.~Svardal, C.~Rueffler, and J.~Hermisson.
\newblock Comparing environmental and genetic variance as adaptive response to
  fluctuating selection.
\newblock {\em Evolution}, 65:2492--2513, 2011.

\end{thebibliography}

\end{document}